\documentclass[twoside,12pt]{article}
\usepackage{amsmath,amssymb,epsfig, amsthm}
\date{}
\newtheorem{Lemma}{LEMMA}[section]
\newtheorem{Corollary}[Lemma]{COROLLARY}
\newtheorem{Theorem}[Lemma]{THEOREM}
\newtheorem{Proposition}[Lemma]{PROPOSITION}
\newtheorem{Definition}[Lemma]{DEFINITION}
\newtheorem{Example}[Lemma]{EXAMPLE}

\newcommand{\bnum}{\begin{enumerate}}
\newcommand{\enum}{\end{enumerate}}
\newcommand{\bi}{\begin{itemize}}
\newcommand{\ei}{\end{itemize}}
\newcommand{\btab}{\begin{tabular}}
\newcommand{\etab}{\end{tabular}}
\newcommand{\beq}{\begin{eqnarray*}}
\newcommand{\eeq}{\end{eqnarray*}}
\newcommand{\beqn}{\begin{eqnarray}}
\newcommand{\eeqn}{\end{eqnarray}}

\newcommand{\bq}{\begin{equation}}
\newcommand{\eq}{\end{equation}}
\newcommand{\CA}{{\cal A}}

\newcommand{\CU}{{\cal U}}

\newcommand{\CH}{{\cal H}}

\newcommand{\CL}{{\cal L}}

\newcommand{\CS}{{\cal S}}
\newcommand{\CT}{{\cal T}}

\newcommand{\CX}{{\cal X}}
\def\phi{\varphi}
\def\epsilon{\varepsilon}
\newcommand{\BS}{\mathbb S}
\newcommand{\BR}{\mathbb R}
\newcommand{\BC}{\mathbb C}
\newcommand{\BF}{\mathbb F}

\newcommand{\Ps}[1]{\mathrm{PG}(#1,\BR)}
\newcommand{\Pd}{\mathrm{PG}(3,\BR)}
\newcommand{\Pf}{\mathrm{PG}(5,\BR)}
\newcommand{\mo}{\mathop\mathrm}
\newcommand{\GL}[1]{{\mo {GL}}_{#1} \BR}

\newcommand{\kasten}{\vbox{\hrule height 8pt width 8.6pt depth -7.4pt
    \hbox{\vrule width 0.6pt height 7.4pt
    \kern 7.4pt \vrule width 0.6pt height 7.4pt}
    \hrule height 0.6pt width 8.6pt}}
\newcommand{\ok}{\hfill\kasten}
\newcommand{\bpf}{\begin{Proof}}
\newcommand{\epf}{\ok\end{Proof}\bigskip\noindent}
\newcommand{\bthm}{\begin{Theorem}}
\newcommand{\ethm}{\end{Theorem}}
\newcommand{\ble}{\begin{Lemma}}
\newcommand{\ele}{\end{Lemma}}
\newcommand{\bprop}{\begin{Proposition}}
\newcommand{\eprop}{\end{Proposition}}
\newcommand{\bcor}{\begin{Corollary}}
\newcommand{\ecor}{\end{Corollary}}

\begin{document}
\title{Regular parallelisms on $\Pd$ from generalized \\ line stars: The oriented case}
\author{Rainer L\"owen}
\maketitle
\thispagestyle{empty}

{\centering\footnotesize \it Dedicated to the memory of Helmut Reiner Salzmann\par}

\begin{abstract}
Using the Klein correspondence, 
regular parallelisms of $\Pd$ have been described by Betten and Riesinger  in terms 
of a dual object, called a hyperflock determining ($\it hfd$) line set. In the special 
case where this set has a span of dimension 3, a second dualization leads to a more 
convenient object, called a generalized star of lines. 
Both constructions have later
been simplified by the author. 

Here we refine our simplified approach in order to obtain similar results for regular 
parallelisms of oriented lines. 
As a consequence, we can demonstrate that 
for oriented parallelisms, as we call them,
there are distinctly more possibilities than in the non-oriented case. The proofs require 
a thorough analysis of orientation in projective spaces (as manifolds and as lattices) 
and in projective planes and, in particular, in translation planes. This is 
used in order to handle continuous families of oriented regular spreads 
in terms of the Klein model of $\Pd$. This turns out to be quite subtle.
Even the definition of suitable classes of dual objects modeling oriented parallelisms
is not so obvious.
\\
 
\end{abstract}

\section{Introduction}\label{intro}
A parallelism in real projective 3-space is an equivalence relation on lines (always
assumed to be continuous in a certain sense) such that every class is a spread 
(i.e., a partition of the point set into disjoint lines). The classical example is 
Clifford parallelism, but there are many more examples with varying amounts of symmetry, as
was shown by Betten and Riesinger; see \cite{Thaswalker}, \cite{hyper} 
and other articles by these authors. 
From every parallelism one can construct 
an oriented parallelism (i.e., a similar equivalence relation on oriented lines), by a 
trivial process of `unfolding'. In \cite{so3}, the author found that 
even among the most symmetric oriented parallelisms other than Clifford parallelism there are 
`non-foldable' examples that do not arise in this trivial way. Our aim is to show that 
a vast number of non-foldable examples can be found among regular oriented parallelisms 
(where the spreads are all isomorphic to the complex spread), even with 2-torus symmetry.  

In the present paper, we lay the foundations for this project by establishing 
powerful construction principles. In the case of non-oriented parallelisms, such 
principles were found by Betten and Riesinger and later simplified by the author. 
One works with the Klein model, which describes $\Pd$ within $\Pf$. The line space of $\Pd$
corresponds to the Klein quadric $K$, and regular spreads become intersections of 
$K$ with certain 3-spaces. 

First, there is a  construction that works for all regular parallelisms; 
by dualizing with respect to the Klein quadric, the parallelism is turned into a set 
$\CH$ of lines not meeting the quadric such that every tangent hyperplane of the 
quadric contains exactly one of them \cite{hyper}, \cite{gldirect}. 
Such sets are called hyperflock determining line sets, abbreviated $\it hfd$ line sets.
In the special case where $\CH$ spans (a plane or) a 3-space, one can dualize once more 
within the 3-space and one obtains a so-called generalized line star or $\it gl$ 
star \cite{Thaswalker}, \cite{gldirect}. 
A $\it gl$ star  is a set of lines in the 3-space 
such that every exterior point with respect to $K$ is 
on exactly one of them. This correspondence has been used to 
construct large sets of examples with torus symmetry \cite{Stevin}, \cite{torus}.

We shall prove similar results for parallelisms of oriented lines 
(briefly called oriented parallelisms). There are several obstacles that have to be overcome. 
This requires a careful analysis of orientation in various geometric contexts, 
see Section \ref{orient}. A spread is a line set homeomorphic to the 2-sphere, 
and the section culminates in the proof that orienting a spread as a manifold 
amounts to the same as orienting all lines in the spread (Theorem \ref{spreadorient}).  
Regular spreads appear in the Klein quadric as intersections with certain 3-spaces. 
In Section \ref{Klein}, we show that this correspondence lifts to a continuous map from 
oriented 3-spaces to oriented spreads, where the topology on the set of oriented spreads is 
defined by a Hausdorff metric. 

In Sections \ref{HFD} and \ref{GL} we define the oriented analogs of $\it hfd$ sets 
and of $\it gl$ stars and obtain our main results. For example, given a 3-space $R$ that meets 
$K$ in an elliptic quadric $Q$, an oriented $\it gl$ star  or $\it gl^+$ star with 
respect to $Q$ is a set 
of oriented secants of $Q$ such that every point $p$ of $R$, not in the interior of $Q$, 
is incident with 
exactly two of them and such that the set of these two oriented lines depends continuously 
on the point $p$. It is the latter condition which makes this approach work. 
Example \ref{patholog} will show that compactness would not do as a surrogate. 
The main results of these two sections 
are Theorems \ref{hfd} and \ref{gl}. They may be summarized as follows.

\bthm\label{SUMMARY}
a) There is a one-to-one correspondence between compact oriented regular parallelisms of 
$\Pd$ and $\it hfd^+$ line sets in $\Pf$.

b) Let $R$ be a 3-space of $\Pf$ intersecting the Klein quadric $K$ in an elliptic 
quadric $Q$. 
There is a one-to-one correspondence between $\it hfd^+$ line sets in $R$ 
and $\it gl^+$ stars with respect to $Q$. 
\ethm

In the final Section \ref{ex} we give criteria
that help to recognize oriented generalized line stars and to construct them, 
and we display non-foldable examples with and without rotational symmetry.
A systematic study of regular oriented parallelisms with 2-torus action will be left 
to a future occasion.\\

\section{Orientation}\label{orient}

We consider various concepts of orientation arising in real projective geometry, and their
relationships. Starting from fairly standard concepts, we proceed to develop 
specialized and not so obvious notions and results concerning orientation in 
projective planes, spreads and parallelisms.
A useful reference to standard notions and constructions related to orientation is \cite{GP}.
All notions, notation, and conventions introduced in this section shall be used 
tacitly in the sequel.

\subsection{Oriented vector spaces and projective spaces}\label{orspac}

As usual, an orientation on the vector space $\BR^k$ is given by an ordered basis
$B = (v_1, ...,y_k)$, and two orientations given by $B_1$ and $B_2$ are considered as 
equal when the linear map sending $B_1$ to $B_2$ has positive determinant. 
A vector space 
$V$ with a fixed orientation will usually be denoted $(V,B)$ or simply $V^+$.
Any basis  defining the same orientation as $B$ will be called a \it positive basis \rm of $(V,B)$.
An oriented differential $k$-manifold is a differential manifold with compatible 
orientations on all its tangent spaces. An orientation in this sense defines also an 
orientation of the underlying topological manifold, i.e., preferred generators for all local 
homology groups in dimension $k$.
An oriented vector space may be considered as an oriented differential  manifold, 
since it coincides with each of its tangent spaces.
\\

The projective space $P_k\BR$ is defined as the quotient space of $\BR^{k+1}\setminus\{0\}$
obtained by identifying a nonzero vector $v$ with every scalar multiple $rv$, $0\ne r \in \BR$.
By an orientation of $P=P_k\BR$ we shall mean an orientation of $\BR^{k+1}$ and denote 
the oriented projective space by $P^+$. 
\\

We stress that this does not mean that we have oriented the differential or 
topological $k$-manifold $P$. 
However, if we restrict the quotient map $\rho: \BR^{k+1}\setminus\{0\} \to P_k\BR$ to 
the unit sphere $\BS_k \subseteq \BR^{k+1}$, we obtain a two-sheeted covering map. The 
nontrivial deck transformation of this covering is the antipodal map $-\rm id$. A given 
orientation of the differential manifold $\BR^{k+1}$ induces an orientation on the unit ball 
(as a manifold with boundary), and this in turn yields an orientation on the boundary 
$\BS_k$ as follows, cp. \cite{GP}.  A basis $B$ of the tangent space $T_s\BS_k$ at a point $s$ is positive if $B$, preceded by an outward pointing vector 
$v \in T_s\BR^{k+1} = \BR^{k+1}$, becomes a positive basis of $\BR^{k+1}$.

Now we try to transfer this orientation to $P_k\BR$ via the map $\rho$ by insisting that 
the tangent map of $\rho$ should preserve orientation on tangent spaces. This works 
without conflict if and only 
if the deck transformation preserves the orientation of $\BS_k$, which is the case if 
and only if $k+1$ is even, i.e., if $k$ is odd. So an odd-dimensional oriented 
projective space as defined here is in fact an oriented manifold, but an 
even-dimensional one is not. \\

\bf Remark. \rm The orientation on the manifold $P_k\BR$ defined above 
does not depend on the choice of a quadratic form defining the unit sphere $\BS_k$.
In other words, if 
$\phi: \BR^{k+1} \to \BR^{k+1}$ is a linear automorphism with positive determinant, then
replacing $\BS_k$ with $\BS_k' = \phi\BS_k$ we get the same result. Indeed, $\phi$ 
is an orientation preserving diffeomorphism between the two spheres, and the quotient 
maps $\rho$ and $\rho'$ 
to the resulting oriented manifolds $P_k\BR$ and $P_k'\BR$ also preserve orientations. 
The induced map $\phi^\flat:  P_k\BR \to P_k'\BR$ satisfies 
$\rho' \circ \phi = \phi^\flat \circ \rho$, and therefore preserves orientations, as well.\\

We note that an orientation of a compact connected 
differential $k$-manifold corresponds to an orientation of the underlying topological 
manifold, which can be represented by a preferred generator of its top homology group,
which is infinite cyclic. This amounts to the same as choosing a preferred generator for each
local homology group in dimension $k$ in a coherent way.

\subsection{Grassmannians and  orientation}

In Incidence Geometry, a projective space is commonly considered as the subspace 
lattice of some vector space with the lattice operations \it join \rm $X\vee Y$ (the span 
of the union of two subspaces $X$ and $Y$) and \it intersection \rm 
$X\wedge Y$. The $k$-dimensional
real projective space, considered as the subspace lattice of $\BR^{k+1}$, will be 
denoted $\Ps{k}$. Its  point set $P_k\BR$ is the set of 1-dimensional vector subspaces.
The projective space $P(X)$ associated with an $(l+1)$-dimensional subspace 
$X \le \BR^{k+1}$ may be considered as a submanifold of $P_k\BR$, homeomorphic to $P_l\BR$. 
Such a subset is considered as an $l$-dimensional projective subspace of $P_k\BR$. 
In this way, $\Ps{k}$ becomes a lattice of subsets of $P_k\BR$. \\

The set of all $l$-dimensional subspaces of $P_k\BR$, or equivalently, the set of 
$(l+1)$-dimensional subspaces of $\BR^{k+1}$, is known as the Grassmann manifold 
${\rm G}_{k+1,l+1}$. Indeed, the transitive action of the general linear group 
$\Delta = \GL{k+1}$
turns this set into a compact, connected differential manifold, namely the 
homogeneous space of the group $\Delta$ modulo the stabilizer $\Delta_X$ 
of any fixed $(l+1)$-dimensional 
subspace $X$. Specifically, this means that the subspace $\delta X$ corresponds to the coset
$\delta \Delta_X$ for every $\delta \in \Delta$. 
The differential manifold ${\rm}G_{k+1,1}$ equals the projective space $P_k\BR$ 
as defined earlier. When we think of the elements of a Grassmann manifold as projective 
subspaces, we prefer to write $P_{k,l}$ rather than ${\rm G}_{k+1,l+1}$. Thus, 
we have $P_k\BR = P_{k,0}$.
For the continuity properties of the lattice operations in $\Ps{k}$, a convenient 
reference is  \cite{Kuehne} or \cite{handb}.\\

Now let us consider oriented subspaces. As before, the set of oriented 
$(l+1)$-dimensional vector subspaces of $\BR^{k+1}$ or, equivalently, the set of oriented 
$l$-dimensional projective subspaces of $P_k\BR$, becomes a compact differential 
manifold when considered as a coset space of $\Delta = \GL{k+1}$. This time, the subgroup 
to be factored out is the stabilizer $\Delta_{X^+}$ of any oriented $(l+1)$-dimensional 
oriented vector subspace $X^+$. The elements of this stabilizer fix $X^+$ as a vector space 
and induce a linear map of positive determinant on $X$. We denote this manifold by
   $$G^+_{k+1,l+1} = P^+_{k,l}.$$
We repeat that it is not possible to consider even-dimensional  projective subspaces as
oriented manifolds. Note also that, in general, we do not have lattice operations for 
oriented subspaces. 

Clearly, the stabilizer $\Delta_{X^+}$ is an index 2 subgroup of $\Delta_X$, and therefore,
the manifold $P^+_{k,l}$ is a 2-sheeted covering space of $P_{k,l}$. 
The covering maps will be denoted $\psi$. We could have used this 
fact as a definition of $P^+_{k,l}$, but this would not allow us to obtain a direct grip on 
the orientation carried by a subspace.

\subsection{Polarities and orientation}

A \it polarity \rm is an antiautomorphism $\xi$ of order 2 of the 
lattice $\Ps{k}$. It is defined 
by a non-degenerate symmetric bilinear form $f$ on $\BR^{k+1}$ and sends a subspace $X$ to
  $$\xi(X) = \{v \in \BR^{k+1} \vert \, f(v,X)=0\}.$$
We have $\dim \xi(X) + \dim X = k+1$ for the vector space dimensions. In terms of projective
dimensions, this means that $\xi$ sends $P_{k,l}$ to $P_{k,k-l-1}$. 

If $f$ restricts to   
a non-degenerate form on $X$, then $\xi(X)$ is indeed a vector space complement of $X$. 
In this case, we can define $\xi$ for oriented subspaces. We choose a fixed orientation for 
$\BR^{k+1}$, and we define  $\xi(X^+)$ by insisting 
that the sum decomposition 
     $$\BR^{k+1} = X \oplus \xi(X)$$
is a sum of oriented vector spaces, i.e., that a positive ordered basis
of $X$, followed by a positive ordered basis of $\xi(X)$, defines the given 
orientation of $\BR^{k+1}$. Thus we have a partial lift of the polarity to the 2-sheeted 
covers, that is, a continuous partial map $\xi^+: P^+_{k,l} \to P^+_{k,k-l-1}$ that 
commutes with the covering maps $\psi$ in the sense that 
$\psi \circ \xi^+ = \xi \circ \psi$.

\subsection{Oriented projective planes}

Let $(P,\CL)$ be a compact, connected topological projective plane. This means first of 
all that $(P,\CL)$ is a projective plane, i.e., that the elements of $\CL$ are subsets of 
the set $P$, called lines, and that two distinct points $p,q$ are joined by a unique line 
$p \vee q \in \CL$ and two distinct lines $K,L$ meet in a unique point $K \wedge L \in P$. 
Moreover, we require that $P$ and $\CL$ are compact, connected topological spaces and that 
the operations $\vee$ and $\wedge$ are continuous. The classical examples are the planes
$\mathrm{PG}(2,\BF)$, where $\BF$ stands for the field of real or complex numbers, the 
skew field of quaternions, or the division algebra of octonions. The examples that we have 
in mind here are translation planes defined by spreads in $\Pd$, see Section \ref{spr}.\\

For simplicity, we shall assume that lines are topological manifolds, which implies that 
they are in fact spheres of dimension $l \in \{1,2,4,8\}$, and that $P$ and $\CL$ are 
$2l$-dimensional manifolds. If we only assume that $P$ 
has finite covering dimension $\dim P < \infty$, then lines are homology $l$-manifolds 
with the same possibilities for $l$ as before, and their homology groups are those of 
an $l$-sphere. Proofs of these facts can be found in Chapter 5 of \cite{CPP}. 
Background information without proofs is also given  in \cite{handb}.\\

Let $K,L$ be two lines and let $x$ be a point not contained in any of these lines.
Then the central projection map
    $$\omega (K,x,L): y \to (y \vee x) \wedge L$$
is a homeomorphism from $K$ to $L$, called a \it perspectivity\rm .  Compositions of 
perspectivities starting from a line $L$ and ending up on the same line are 
called \it projectivities\rm . They form a group $\Omega(L)$. These groups have been 
studied extensively \cite{wind}.\\

Every line $L \approx \BS_l$ admits two possible orientations, and we denote the set of 
all oriented lines by $\CL^+$. There is a two to one surjective map 
    $$\psi: \CL^+ \to \CL$$
that forgets orientations. Orientation forgetting maps will occur frequently, and they will 
always be denoted $\psi$. Our goal is to define a topology on $\CL^+$ such that $\psi$ becomes 
a two-sheeted covering map. The first such construction was given by Salzmann \cite{adv}, p. 10,
using the space of all embeddings $\BS_l \to P$ whose images are lines, 
equipped with the compact open topology. Here we use a different approach.
\\

Let $L^+$ be an oriented line, and choose a point $x \notin L = \psi(L^+)$. 
Define a set of oriented lines
   $$\CL^+(L^+,x)$$
to be the set of all lines $M$ not containing $x$, endowed with the orientation 
transferred from $L^+$ via $\omega(L,x,M)$.  The restriction of $\psi$ to this set is a 
bijection onto an open subset of $\CL$ (the set of lines not containing $x$), and we 
define a topology on $\CL^+(L^+,x)$ by insisting that this bijection be a homeomorphism.

Now consider two such sets $\CU_1 = \CL^+(L_1^+,x_1)$ and $\CU_2=\CL^+(L_2^+,x_2)$, and 
let $\CX^+$ be the set
of oriented lines containing neither $x_1$ nor $x_2$. An oriented line
$M^+ \in \CU_1$ belongs to the intersection $\CU_1\cap \CU_2$ if and only if it 
lies in $\CX^+$ and 
the map $\omega(L_1,x_1,M,x_2,L_2)$ (to be read as a composition of 
perspectivities in the obvious manner) is orientation preserving as a map $L_1^+ \to L_2^+$. 
If $M_t$ is a path in $\CX = \psi\CX^+$, then the corresponding maps $\omega_t$ 
are homotopic, hence they have the same effect on the top homology groups and
are either all orientation preserving or all orientation reversing.
Thus the intersection $\CU_1\cap \CU_2$ is mapped by the forgetful 
map $\psi$ onto a (possibly empty) union of some path connceted components of $\CX$. 
These components are 
open sets, so we see that $\CU_1 \cap \CU_2$ is open in $\CU_1$ and in $\CU_2$ and inherits the 
same topology from both sets. 

Now we endow $\CL^+$ with the topology generated by all the topologies 
on the various sets $\CL^+(L^+,x)$, and it is still true that $\psi$ restricts to a 
homeomorphism on each of these sets with respect to this topology. Thus we have obtained

\bprop
If the set $\CL^+$ of oriented lines of a compact projective plane is equipped with the 
topology defined above, then the forgetful map $\psi: \CL^+ \to \CL$ becomes a 
two-sheeted covering map. \ok
\eprop

From this, we obtain the next proposition almost as a corollary. By a \it section \rm of 
the map $\psi$ we mean a map $\sigma$ in the opposite direction such that the composition 
$\psi \circ \sigma$ is the identity map of $\CL$.

\bprop
\item{a)} Let $(P,\CL)$ be a projective plane with lines of dimension $l$.
The space $\CL^+$ of oriented lines is connected if $l = 1$ (in fact, 
it is then a 2-sphere) and is disconnected otherwise.
\item{b)} The forgetful map $\psi: \CL^+ \to \CL$ admits a section if and only if $l \ge 2$.
\eprop

\bpf
If $l \ge 2$, then the point space $P$ is simply connected, because the complement of a 
point $P\setminus \{x\}$ deformation retracts onto every line not containing $x$, 
cp. \cite{CPP}, 51.26.
Exchanging the roles of points and lines, we see that $\CL$ is simply connected as well. 
Therefore, a two-sheeted covering of $\CL$ must be the topological sum of two copies of $\CL$, each 
of which is mapped homeomorphically onto $\CL$ by the covering map.\\

For $l=1$, the space $\CL^+$ is connected. Indeed, the pencil $\CL^+_x$ of oriented 
lines passing through a given point $x$ is connected (in fact, homeomorphic to $\BS_1$) 
since we have a continuous surjection from the boundary of a disc containing $x$ in 
its interior 
onto $\CL_x^+$, sending a point $p$ to the line $p \vee x$ oriented locally from $p$ to $x$. 
Furthermore, any two oriented lines belong to the pencil of oriented lines 
passing through their 
point of intersection.

Now the space $\CL$ is homeomorphic to $P_2\BR \approx \BS_2/\pm{\rm id}$ for $l=1$, 
see \cite{CPP}, 42.10. 
Thus, the 2-sphere is the only connected two-sheeted covering space of $\CL$. 
A section  to $\psi$ would be a homeomorphism by domain invariance, a contradiction.
\epf

For our purposes, the preceding result is not enough. We need an explicit construction 
of sections in the case $l \ge 2$. This will be made possible with the aid of 
projectivity groups $\Omega(L)$. \\

It is rather easy to see that projectivities do not preserve
orientation of lines if $l=1$. For example, in the real affine plane (which canonically extends
to the projective plane) consider the lines $X$ (the $x$-axis) and $Y$ (the $y$-axis) and 
the points $p = (1,-1)$ and $q =(-1,-1)$. The projectivity $\omega(X,p,Y,q,X)$ (to be read as 
a composition of perspectivities in the obvious manner) reverses the orientation of the 
$x$-axis. However, for $l \ge 2$, we shall prove that all lines can be oriented in such a 
way that all perspectivities and, hence, all projectivities preserve these orientations. 
By the definition of the topology on $\CL^+$, this then provides a section to the 
forgetful map $\psi$.\\

As a preparation, we endow every group $\Omega(L)$ of projectivities with the compact open 
topology, which turns it into a topological transformation group of the space $L$. 
Next we note that a path $\sigma: [0,1] \to \Omega(L)$ corresponds to an isotopy on $L$,
and so the maps $\sigma(0)$ and $\sigma(1)$ are either both orientation preserving or 
both orientation reversing, because they have the same effect on the top homology group. 
Together with the above example, this shows that $\Omega(L)$ is not pathwise connected in 
the case of the real projective plane. However, we have the following
well-known

\ble
For $l \ge 2$, all groups $\Omega(L)$ of projectivities are pathwise connected.
\ele

\bpf
We follow \cite{wind}, Theorem 3.3. Consider a projectivity 
     $$\omega = \omega(L_1,x_1,L_2,x_2,...,L_{n-1},x_{n-1},L_n)$$ 
with $L_n = L_1=L$. We shall construct a path in $\Omega(L)$ joining 
$\omega$ to the identity. We choose a point $x$ not on any of the lines $L_i$ and join 
every $x_i$ to $x$
by a path $x_i(t)$ in the connected set $(x \vee x_i) \setminus (L_i \cup L_{i+1})$ 
(or by a constant path if $x = x_i$). Replacing every $x_i$ by $x_i(t)$ in the definition 
of $\omega$ we obtain 
a path $\omega(t)$ in $\Omega(L)$. Since all projection centers of $\omega(1)$ are 
identical, $\omega(1)$ is the identity, and $\omega(0) = \omega$.
\epf

Now we obtain a geometric construction of sections of the forgetful map $\psi$.

\bthm\label{planeorient}
In a compact projective plane of dimension $2l \ge 4$ it is possible in exactly two ways to 
orient all lines in such a way that all perspectivities are orientation preserving.
\ethm

\bpf
Start with one line $L$ and orient it in one of the two possible ways. The condition 
on perspectivities then forces our way of orienting all other lines. If conflicts  
arise when we 
transfer orientations about, then it means that two projectivities from $L$ to some 
line $K$ take the orientation of $L$ to distinct orientations of $K$. Then the quotient 
of these projectivities
is an orientation reversing projectivity of $K$ to itself. This is impossible, because
$\Omega(K)$ is path connected, which implies that every projectivity of $K$ to itself 
is isotopic to the identity and hence induces the identity on the top homology group. 
\epf

Note that the orientations constructed in this proof depend continuously on the lines, by 
the very definition of the topology on $\CL^+$. Thus we have in fact obtained a 
section to the forgetful map.

\subsection{Oriented spreads}\label{spr} 

We are now ready to study the orientation properties of spreads of $\Pd$, which are 
crucial to us 
since  parallelisms are built from spreads. In fact, the results of this section, 
together with the continuity result Theorem \ref{orsprbyinters} below, constitute 
the most subtle steps 
towards our final goals. 

When we first introduced oriented parallelisms 
in \cite{so3}, we were content with a very simple approach. A spread is a certain set $\CS$ of 
lines of a projective 3-space, which is homeomorphic to the 2-sphere. Hence, the two-sheeted 
covering $\CS^+$ is necessarily disconnected, and we defined an orientation of $\CS$ to be a 
choice of one of the two connected components of this cover. In the present situation, we need  
closer control of the orientations of the lines $L\in \CS$ involved here, so we need to refine our definition. \\

First we recall the definition of a spread; compare Section 64 of \cite{CPP}. 
Consider the line space $P_{3,1}$ of $\Pd$. A compact set 
$\CS \subseteq P_{3,1}$ is a (topological) \it spread \rm if each point $x$ belongs to a unique 
element  $S_x \in \CS$. By compactness, the map $x \to S_x$ is continuous, and $\CS$ is homeomorphic to 
the 2-sphere (this will become apparent later). \\

One reason for studying spreads is that a spread $\CS$ defines an \it affine translation 
plane \rm 
$\CA_\CS$. We may consider $\CS$ as a subset of the Grassmann manifold $G_{4,2}$, 
i.e., as a set of 2-dimensional vector subspaces of $\BR^4$. The plane $\CA_\CS$ has 
point set $\BR^4$, and its lines are all translates
$L+v$, where $L\in \CS$ and $v \in \BR^4$. Thus, $\CS$ is the pencil of lines 
containing the origin $0\in \BR^4$ (which incidentally explains why $\CS \approx \BS_2$). \\

The projective plane associated with $\CA_\CS$ is the \it projective translation plane \rm 
$\CT_\CS$ defined by $\CS$. There is a particularly nice description of this plane, which is somewhat hidden in the proof of Theorem 64.4 of \cite{CPP}. (That book chapter contains a 
thorough introduction to translation planes.) The construction of $\CT_\CS$ starts from $\Ps{4}$. 
Choose a hyperplane $H$ of that projective space. Then $H$ is isomorphic to $\Pd$, 
and we may consider $\CS$ as a set of lines of $H$. We simply 
write $P_4$ for the point set $P_{4,0}=P_4\BR$ of $\Ps{4}$. 
Now the points of $\CT_\CS$ are
\begin{itemize}
\item the points of $P_4$ not belonging to $H$ and
\item the elements of $\CS$. 
\end{itemize}
This partitions the point set $P_4$ into disjoint sets 
(most of them singletons), and the topology of the point set $T$ of $\CT_\CS$ is 
the resulting quotient topology. 

The lines of $\CT_\CS$ are 
\begin{itemize}
\item the elements of $P_{4,2}$ (2-spaces in $\Ps{4}$) that meet $H$ in a line 
$S \in \CS$ and 
\item the set $\CS$ itself. 
\end{itemize}
Incidence of points and lines is given by inclusion, and 
the topology of the line set is again a suitable quotient topology. \\

Now it is important to note that the topology of $\CS$ considered as a set of 
points of $\CT_\CS$ 
is the same as the topology of $\CS$ inherited from $P_{3,1} = G_{4,2}$, 
the line space of $H$. By definition, $\CS$ 
considered as 
the point set of the line $\CS$ of $\CT_\CS$ is the quotient of the point set 
$H_0$ of $H$ with respect to the map $H_0 \to \CS$ that sends a point to the 
unique spread line containing it. As we noted earlier, this map is continuous 
when $\CS$ is given the topology coming from the line 
space of $H$. It is also closed (by compactness) and surjective, hence it is a 
quotient map and our claim follows. \\

Remembering Theorem \ref{planeorient}, we now conclude that orienting a spread $\CS$
(the special line of $\CT_\CS$) as a manifold amounts to the same thing as orienting \it all \rm
\, lines of $\CT_\CS$. The affine plane $\CA_\CS$ is obtained from $\CT_\CS$ by deleting 
the special line $\CS$ and all its points, and so by orienting all lines of $\CT_\CS$ 
we have in particular oriented all elements of $\CS$, considered as 2-spaces in 
$\BR^4$ (because they are lines of $\CA_\CS$). \\

In terms of the affine plane, this transfer of orientations is easy to visualize: 
$\CS$ is the pencil of lines passing through the origin, and in order to orient $S\in \CS$, 
apply to $S$ a translation $s \to s+v$ with $v \notin S$. Then consider the bijection
    $$s+v \to (s+v)\vee 0$$
of $S+v$ onto $\CS \setminus \{S\}$ and transfer the orientation via these two maps.  
This is, indeed, much simpler, but in order to prove consistency by applying
Theorem \ref{planeorient}, we prefer the projective version.  \\

Summarizing these constructions, we obtain the following theorem.
 
\bthm\label{spreadorient}
Let $\CS \subseteq P_{3,1}$ be a compact spread of $\Pd$. The construction given above 
defines a bijective correspondence between orientations 
of the manifold $\CS \approx \BS_2$ on the one hand and coherent orientations of all lines of 
the projective translation plane $\CT_\CS$ associated with $\CS$ on the other hand, i.e., 
sections to the orientation-forgetting map $\psi$ of the latter plane. 

In particular, orienting $\CS$ as a manifold amounts to the same as orienting, in a coherent way, all the vector 2-spaces $S \in \CS$ as manifolds (or as vector spaces). \ok
\ethm

\subsection{Oriented parallelisms}

A \it parallelism \rm $\Pi$ on $\Pd$ is a set of pairwise disjoint spreads 
covering the line set 
$P_{3,1}$. In other words, a parallelism partitions the line set and therefore
is often thought of as an equivalence relation on the line set. Likewise, 
a \it parallelism of oriented lines \rm or briefly, an \it
oriented parallelism \rm $\Pi^+$ is defined as a set of oriented spreads partitioning the 
set $P^+_{3,1}$ of oriented lines. If $\Pi$ is a parallelism, then an oriented 
parallelism $\Pi^+$ can be obtained from it by taking all oriented spreads $\CS^+$ such that 
$\psi \CS \in \Pi$, where as always $\psi$ denotes the map that forgets orientations. 
This process will be called \it unfolding, \rm and 
oriented parallelisms obtained in this way will be said to be \it foldable\rm. 
Our investigation 
is motivated by the existence of non-foldable oriented parallelisms with nice properties, 
the discovery of 
which is described in \cite{so3}, and by the desire to find more non-foldable 
examples worth looking at. 

In a topology to be introduced shortly, an ordinary parallelism is a 
non-orientable 2-manifold (a projective plane) and an oriented parallelism is a 
2-sphere.
To avoid possible confusion, we stress here that 
the possible orientations of this sphere are irrelevant to us. 
So the term `oriented parallelism' is merely a shorthand for `parallelism of oriented lines'.
This contrasts with the situation for spreads, compare the preceding subsection. 

For the convenience of the reader, we recall 
some basic facts obtained in \cite{so3}, adding some details that require extra 
attention in the present situation. \\ 

We need a condition to ensure topological well-behavedness of a parallelism $\Pi$ or $\Pi^+$. 
A good choice for a topology on $\Pi$ or $\Pi^+$
is the topology  defined by the \it Hausdorff metric \rm \cite{Tuzh},
which we introduce next. 
Let $(X,d)$ be a compact metric space. The \it hyperspace \rm $h(X)$ is the set of 
all compact subsets of $X$, endowed with the metric
   $$d_h(A,B) = \max\{\max_{a\in A}d(a,B), \max_{b \in B}d(b,A)\},$$
where, as usual, $d(a,B) = \max_{b \in B} d(a,b)$. By the \it Hausdorff topology, \rm
we mean the topologgy on $h(X)$ induced by the Hausdorff metric.
Two metric topologies agree if they 
induce the same notion of convergence. In the case of the Hausdorff metric, this notion 
is captured by the following Lemma, which follows easily from the definition, 
in view of compactness. 

\ble\label{hausdconv}
Let $(X,d)$ be a compact metric space. A sequence $A_n \in h(X)$ converges to $A \in h(X)$ 
if and only if the following two conditions are satisfied.
\item{i)}  If a sequence of points $a_n \in A_n$ converges to a point $a$ in $X$, then 
$a \in A$. 
\item{ii)} Every point $a \in A$ is the limit of some sequence of points $a_n \in A_n$.
\ok
\ele 

In what follows, we want to treat ordinary and oriented parallelisms simultaneously. We write
$\Pi^*$ and $L^* \in P_{3,1}^*$ to indicate that we are thinking of both possibilities. By 
      $$\Pi^*(L^*)$$ 
we denote the unique spread in $\Pi^*$ that contains $L^*$. Similarly, 
    $$\Pi^*(x,L^*)$$ 
denotes the unique line that belongs to the same spread as $L^*$ and contains a given point $x$.
In this way, we define  two maps 
$P_{3,1}^* \to \Pi^*$ and $P_{3,0}\times P_{3,1}^* \to P_{3,1}^*$; each of them 
contains the same information as $\Pi^*$ itself, which should justify the abuse 
of notation. We have the following

\bprop\label{toppar}
Let $\Pi^*$ be an ordinary or oriented parallelism on $\Pd$. The following conditions 
are equivalent.
\item{1)} $\Pi^*$ is compact with respect to the Hausdorff topology on $h(P_{3,1}^*)$.
\item{2)} The map $\Pi^*: P_{3,1}^* \to h(P_{3,1}^*)$  defined above is continuous with respect to the Hausdorff topology on the hyperspace.
\item{3)} The map $\Pi^*: P_{3,0}\times P_{3,1}^* \to P_{3,1}^*$ defined above is continuous. 
\eprop

\bpf
Assume (1). In order to show (3), we prove sequential continuity. 
If $(x_n,L_n^*) \to (x,L^*)$, we have to show that $\Pi^*(x_n,L_n^*) \to \Pi^*(x,L^*)$. By compactness of $\Pi^*$, we may assume that 
$\Pi^*(L_n^*)$ converges to some spread $\CS^* \in \Pi^*$. 
Then by Lemma \ref{hausdconv}, we have
$L^* \in \CS^*$, and so $\CS^* = \Pi^*(L^*)$. Moreover, we may assume that $\Pi^*(x_n, L_n^*)$
converges to some line $K^* \in P^*_{3,1}$.  Then $x \in K^*$ since incidence is closed, 
and $K^* \in \Pi^*(L^*)$ by Lemma \ref{hausdconv}, because $\Pi^*(x_n,L_n^*) \in \Pi^*(L_n^*)$
and $\Pi^*(L_n^*) \to \Pi^*(L^*)$. Thus, $K^* = \Pi^*(x,L^*)$. This proves (3).

Now assume (3) and suppose that $L_n^* \to L^*$. In order to prove (2), 
we have to show that $\Pi^*(L_n^*)\to
\Pi^*(L)$. Thus we have to verify conditions (i) and (ii) of Lemma \ref{hausdconv}. 
So let $K^*_n\in \Pi^*(L_n^*)$ and assume that $K_n^* \to K^* \in P^*_{3,1}$. 
We have to show that $K^* \in \Pi^*(L)$. Choose points $x_n \in K_n^*$. We may assume that 
$x_n \to x \in P_{3,0}$. Then by (3), we have
   $$K_n = \Pi^*(x_n,L_n^*) \to \Pi^*(x,L^*) \in \Pi^*(L^*).$$
This proves (i). For condition (ii), let $K^* \in \Pi^*(L^*)$. We are looking for lines
$K_n^* \in \Pi^*(L_n^*)$ such that $K_n^* \to K^*$. Choose any point $x \in K^*$. Then 
$K_n := \Pi^*(x,L_n^*)$ belongs to $\Pi^*(L_n^*)$, and these lines converge to 
$\Pi^*(x,L^*) = K^*$ by (3).

Finally, (2) implies (1) because the map considered in (2) restricts to a bijective map 
from the star of all lines or oriented lines containing any chosen point $x$ to $\Pi^*$. The star is compact, and (1)
follows.
\epf

We shall say that $\Pi^*$ is a \it topological (oriented) parallelism \rm if it satisfies these 
equivalent conditions.
Line stars are homeomorphic to $P_2\BR$ in the ordinary case and to $\BS_2$ in the oriented case. Hence, the last step of the above proof shows: 

\ble
A topological parallelism $\Pi$ is 
homeomorphic to $P_2\BR$, and a topological oriented parallelism $\Pi^+$ is homeomorphic to $\BS_2$.
\ok
\ele

\section{Klein correspondence for oriented regular spreads}\label{Klein}

Within $\Ps{5}$, the Klein correspondence sets up a model of $\Pd$ that is well suited for studying the line space $P_{3,1}$. We summarize without proofs properties of the Klein correspondence that can be found
in the literature, in particular in \cite{Knarr}, \cite{Stevin}, \cite{gldirect}. 
The Klein model arises from the index 3 bilinear form 
$f(x,y) = x_1y_1 + x_2y_2 +x_3y_3 - x_4y_4-x_5y_5-x_6y_6$ on $\BR^6$. There are several 
kinds of special subspaces of $\BR^6$ with respect to this form:
\begin{itemize}
\item Totally isotropic one-dimensional subspaces. Viewed as points of $\Ps{5}$, they constitute 
the \it Klein quadric \rm $K$, which represents the line set of $\Pd$.
\item Two sorts of totally isotropic 3-dimensional subspaces (i.e., $f$ induces the zero 
form on them). They represent the points and hyperplanes of $\Pd$. Incidence with lines is represented as reverse or direct inclusion, respectively. 
\item Four-dimensional subspaces of signature $(1,3)$ or $(3,1)$, i.e., $f$ induces a non-degene\-rate form of index one on them. Viewed as projective subspaces, they are 
3-dimensional and intersect $K$ in an elliptic quadric. They will be most important to us, 
and we call them (1,3)-spaces or (3,1)-spaces, respectively. As a vector space, a (1,3)-space 
contains a 3-dimensional negative definite subspace, and a (3,1)-space contains a 3-dimensional
positive definite subspace.
\end{itemize}

We use the Klein correspondence mainly to describe regular spreads. A \it regular 
spread \rm of $\Pd$ is a spread isomorphic to the \it complex spread\rm , which defines the complex affine plane. In other words, the latter spread consists of the one-dimensional complex 
subspaces of $\BC^2 = \BR^4$. We have the following Lemma, proved, e.g., 
in \cite{Stevin}, Proposition 13.

\ble\label{regspr}
In the Klein model of $\Pd$, the regular spreads are precisely the sets $\CS=K\cap P$, 
where $P$ is a $(1,3)$-space or a $(3,1)$-space.
\ele

Here, and frequently in what follows, we identify projective subspaces 
with their point sets, so that the lattice $\Pf$ is viewed as a lattice of subsets of $P_5\BR$.
Our aim now is to obtain an oriented version of the above lemma, together with a 
continuity assertion. This will be achieved by the  construction given below.\\

The space of oriented lines of $\Pd$ is a two-sheeted covering of the line space, 
which we identify with the Klein quadric $K$. Therefore, 
we shall write $K^+$ for the space of oriented lines
and call it the \it oriented Klein quadric. \rm 
The covering map $K^+ \to K$ will be denoted $\psi$, as usual. It is well-known that 
$K^+ \approx \BS_2 \times \BS_2$. An easy proof can be given using the oriented left and right
Clifford parallelisms, see \cite{so3}, Proposition 2.4. \\

\bf Construction.  \it Step 1.  Let $P^+ \in P^+_{5,3}$ be an oriented (3,1)-space, 
considered as projective space. The projective dimension of $P^+$ is odd, so, 
as explained in Section \ref{orspac}, we are given an orientation of 
$P^+$ as a differential manifold. 

\it Step 2.  Let $P = \psi(P^+)$. The elliptic quadric $\CS = P \cap K$, 
which represents a regular spread,
separates $P$ into two components, and is the boundary of both. 
The closure of the `interior' component is a 
compact ball and inherits an orientation from $P^+$. This orientation induces an 
orientation on the boundary $\CS$ as follows: Let  $B=(v_2,v_3)$ be an ordered 
basis of the tangent space  $T_s\CS$ at $s \in \CS$ and let $v \in T_sP$ 
be an outward pointing tangent vector.
Then the orientation of $\CS$ at $s$ is defined by $B$ if $(v,v_2,v_3)$ 
is a positive basis of $T_sP$. 

\it Step 3.  Finally, we apply Theorem \ref{spreadorient}, and from the orientation of 
the manifold $\CS$ we obtain an orientation of all lines of the spread, such that we end up 
with one particular connected component $\CS^+$ of $\psi^{-1}(\CS) \subseteq K^+$. 
We shall denote this set by
      $$\CS^+ = \CS^+(P^+) \subseteq K^+.$$

\rm Our aim is to show now that the map $P^+ \to \CS^+(P^+)$ from $P_{5,3}^+$ to the 
hyperspace $h(K^+)$
is continuous. This is the main step in the proof of Theorem \ref{hfd} below, which 
describes a way of constructing all oriented regular parallelisms. 
We need some preparation concerning group actions. For more details on the group actions 
discussed here, see \cite{Stevin}. \\

We consider the group
   $$\Sigma = \mathop{\rm {SL}}(4,\BR).$$ 
By definition, this group acts on $\BR^4$; it also acts (ineffectively) on $\Pd$. 
Via the Klein correspondence, the latter action is translated to an $f$-orthogonal 
action on $\BR^6$ and on 
$\Ps{5}$, where, as earlier, $f$ denotes the form defining the Klein quadric. 
In fact, $\Sigma$ induces an index 2 subgroup of the projective orthogonal group.
In particular,
the action of $\Sigma$ leaves invariant all the sets of subspaces  of special type enumerated
at the beginning of this section. The actions on these sets of subspaces are transitive. 
The action also lifts to the sets $P^+_{5,l}$ of oriented 
subspaces.  \\

A transitive action of a Lie group $\Gamma$ on a manifold $M$ always 
admits \it local sections\rm.
That is, given a point $x \in M$, there exist a neighborhood $U$ of $x$ and a continuous 
map $u \to \gamma_u$ from $U$ into $\Gamma$ such that $\gamma_x$ is the 
identity and $u= \gamma_u(x)$ for all $u\in U$. 
This can be shown by proving that the map $\Gamma \to M$ sending $\gamma$ to $\gamma(x)$ 
is a submersion. See \cite{naga} for a more general result.

\ble\label{section} 
The actions of $\Sigma$ on the sets of oriented and non-oriented special 
subspaces of $\Ps{k}$ listed at the beginning of this section admit local 
sections. \ok
\ele

\ble \label{hausdconv2}
Let $\Gamma$ be a topological group acting on a metric space $X$ and let $Y \subseteq X$ 
be compact. If $\gamma_\nu$ is a sequence in $\Gamma$ converging to the unit element, then 
$\gamma_\nu(Y) \to Y$ in the Hausdorff metric.
\ele

\bpf The continuity of the map $\Gamma\times X \to X$ 
sending $(\gamma,x)$ to $\gamma(x)$ implies that 
$\gamma_\nu \to {\rm id}$ uniformly on the compact set $Y$.  The claim follows easily.
\epf

\ble\label{grassm}
If $\Ps{k}$ is viewed as a lattice of subsets of $P_k\BR$, 
then the topology of the Grassmannians
$P_{k,l}$ is induced by the Hausdorff metric. 
\ele

\bpf
We use the continuity properties of the topological projective space $\Ps{k}$ given, 
e.g., in \cite {Kuehne} or \cite{handb}.
Suppose that $P_\nu \to P$ in $P_{k,l}$. Choose a subspace $Q$ of dimension $k-l-1$ 
that is in general position to all these spaces. Then the projection map $\pi_\nu: P \to
P_\nu$ sending $p \in P$ to $(p\vee Q)\wedge P_\nu$ is a homeomorphism, and 
$\pi_\nu$ converges to the identity of $P$. Hence for every point $p\in P$, the sequence
$\pi_\nu p_\nu \in P_\nu$ converges to $p$, and condition (ii) of Lemma \ref{hausdconv} 
is satisfied.  Condition (i) follows from the fact that the incidence relation is closed.
This shows that $P_\nu \to P$ with respect to the Hausdorff metric. Alternatively, this can be 
deduced from the preceding two lemmas.

Conversely, assume that $P_\nu \to P$ in the Hausdorff sense. Let $X\subseteq P$ 
be a set of $l+2$ points spanning $P$ (a \it frame \rm for $P$). Then for every $\nu$ 
there is a set $X_\nu \subseteq P_\nu$ of cardinality $l+2$ such that $X_\nu \to X$ in 
the Hausdorff sense. It follows that for $\nu$ large enough, $X_\nu$ is a frame for $P_\nu$,
and the continuity of forming spans implies that $P_\nu \to P$ in $P_{k,l}$.
\epf

\bthm\label{orsprbyinters}
The above construction defines a continuous map $P^+ \to \CS^+(P^+)$ from the set of 
oriented $(3,1)$-spaces in $P^+_{5,3}$ to the hyperspace $h(K^+)$ of the 
oriented Klein quadric. 
\ethm

\bpf
1. Let $P_\nu^+ \to P^+$ be a convergent sequence of $(3,1)$-spaces in $P^+_{5,3}$. 
According to Lemma \ref{section}, 
there is a sequence $\sigma_\nu \in \Sigma$, converging to the identity, such that 
$P^+ = \sigma_\nu P^+_\nu$ for all $\nu$. Using Lemma \ref{hausdconv2}, we infer that that 
$\psi P^+_\nu = P_\nu$
converges to $\psi P^+ = P$ in the Hausdorff metric; compare also 
Lemma \ref{grassm}. Since $K$ is invariant 
under $\Sigma$, the same also holds for
the intersections with $K$, that is, $\CS_\nu = P_\nu \cap K \to \CS = P\cap K$ 
in the hyperspace $h(K)$. 

2. The group $\Sigma$ of isomorphisms respects all steps of the construction that turns 
an oriented projective subspace of odd dimension into an oriented manifold. Therefore, the oriented manifolds $P_\nu^+$ converge to the oriented manifold $P^+$. 
By this we mean that  that there is a 
uniformly continuous sequence of orientation preserving 
embeddings of $P_3^+\BR$ onto 
$P_\nu^+\subseteq P_5$ that converges to an orientation preserving embedding onto 
$P^+$. Indeed, the restrictions of the maps $\sigma_\nu^{-1}$ form such a sequence.

3. The group $\Sigma$ respects all structural features that were used in the construction of an
orientation of the spreads $\CS_\nu$ and $\CS$ from the orientations of $P_\nu$ and $P$,
given in Step 2 of the construction following Lemma \ref{regspr}. As before, this 
implies that $\CS_\nu^+ \to \CS^+$ as oriented manifolds.  

4. The group $\Sigma$  acts on both $\BR^4$ and  $\Pd$.
The group elements  $\sigma_\nu$ send the spread $\CS_\nu$, considered as a
subset of $G_{4,2}$, to the spread $\CS$. Consequently, $\sigma_\nu$ is an isomorphism between 
the affine translation planes defined by these spreads, and extends to an 
isomorphism of the associated projective planes. Thus $\sigma_\nu$ preserves all 
steps in the construction of orientations on the elements of $\CS_\nu$ and 
of $\CS$. Consequently, $\sigma_\nu$ 
sends $\CS^+(P_\nu^+)$ to $\CS^+(P^+)$. Moreover, $\sigma_\nu$ 
converges to the identity on $K^+$, 
hence as before we may conclude that $\CS^+(P_\nu^+) \to \CS^+(P^+)$ in the 
Hausdorff metric, as desired.
\epf

\section{Oriented $\it hfd$ line sets and first main result}\label{HFD}

We begin with a topological Lemma. It is probably known, but I do not remember seeing it in the literature. If $X$ is a topological space, we let $h_2(X)$ denote the set of subsets $A \subseteq X$
with cardinality $\# A = 2$. We topologize this as the quotient
  $$h_2(X) = \left ( (X\times X) \setminus \Delta_X \right ) /\langle s \rangle,$$
where $\Delta_X$ denotes the diagonal $\{(x.x) \vert \, x \in X\}$ and $\langle s \rangle$ 
is the the group generated by the switching map
$s: (x,y) \to (y,x)$. If $X$ is metric, then this topology is also induced by the Hausdorff metric, which is why we choose the symbol $h_2$. 

\ble\label{2-1}
Let $q: \tilde Y \to Y$ be a two-sheeted covering map of connected  Hausdorff spaces, and let $X$ be a compact space. Let $\tilde g:X \to \tilde Y$ be continuous and $g = q\circ \tilde q$.

Suppose that all inverse images $g^{-1}(y)$, $y\in Y$, have cardinality 2 and that the resulting map
$g^{-1}:Y \to h_2(X)$ is continuous. 

Then the map $\tilde g$ is bijective and, in fact, a homeomorphism.
\ele

\bpf
Let $U\subseteq Y$ be an open set which is evenly covered, that is, $q^{-1}(U)$ 
is a union of two open subsets $U_1, U_2$ which are both mapped homeomorphically 
onto $U$ by $q$. 
If $\tilde g$ maps the two $g$-inverse images of $u\in U$  into the same sheet $U_i$, 
then the images are in fact equal, and
the same happens for nearby points $u'$, by  continuity of $g^{-1}$. 
On the other hand, if those images are distinct, then the same holds in a neighborhood of $u$. 
This shows that the cardinality of $\tilde g (g^{-1}(y))$, $y \in Y$, is locally constant. 
By connectedness of $Y$, this cardinality is always 1 or always 2. In the latter case, 
$\tilde g$ is bijective and hence a homeomorphism by compactness. In the former case, 
looking at the sheets again one sees that
the set $\tilde g(X)$ and its complement are both open and nonempty, 
a contradiction to connectedness.  
\epf

\begin{Example}
\rm The following shows that the assumption about continuity of the inverse 
is indispensable in the Lemma above. We take $\tilde Y = X = \BS_2$ and $Y= P_2\BR$, the real projective plane. There is the two-sheeted covering $q: \tilde Y \to Y$ 
sending $y$ to $\pm y$. Consider $g = q\circ \tilde g$, 
where $\tilde g$ is either the identity map or the folding map 
$(x,y,z) \to (x,y,\vert z \vert)$. In the second case, $g^{-1}$ is discontinuous at the 
equator ($z=0$), and $\tilde g$ is neither injective nor surjective.
\end{Example}

We now return to the Klein model of $\Pd$, and consider the polarity $\pi_5$ defined by the 
bilinear form $f$ of signature $(3,3)$. First we note that the polar $\pi_5(x)\in P_{5,4}$ 
of a point $x \in K$ (which represents a line of $\Pd$) is the \it tangent hyperplane \rm of 
$K$ at that point. This is not to be confused with the tangent vector space $T_xK$ of the differential manifold $K$,
which is of no concern to us at the moment. 

A line $L \in P_{5,1}$ is called an \it exterior line \rm  with respect to 
$K$ if $L\cap K = \emptyset$. We note that this is the case if and only if $L$ is a 
subspace  of type $(2,0)$ or $(0,2)$, i.e., the form $f$ is positive or negative 
definite on $L$ considered as a two-dimensional vector space. 
Then $\pi_5(L) \in P_{5,3}$ is a $(1,3)$-space or a $(3,1)$-space, respectively, and defines a regular spread $\pi_5(L)\cap K$.  \\

Betten and Riesinger \cite{hyper}
defined a \it hyperflock determining line set \rm or shortly, an \it hfd set 
\rm to be a set $\CH \subseteq P_{5,1}$ of exterior lines such that every tangent hyperplane
$\pi_5(x)$, $x \in K$, contains exactly one line from $\CH$. Since $\pi_5$ 
is an antiautomorphism
of the lattice $\Ps{5}$ (i.e., it reverses inclusions), this implies that every 
$x\in K$ is contained in exactly one element of $\pi_5(\CH)$. In other words, this set  
defines a regular parallelism $\Pi (\CH)$ by taking intersections with $K$. 
Moreover, every regular 
parallelism  arises in this way, and the parallelism $\Pi (\CH)$ is topological 
if and only if $\CH$ is compact; compare \cite{gldirect}. Now we imitate this 
in the oriented case, but we need to change in the pattern. For example, 
inclusion is only defined for non-oriented subspaces, so we cannot directly mimick 
the above definition of an \it hfd \rm line set. 

\begin{Definition}
\rm a) An \it  oriented hfd line set \rm or briefly, an  $\it hfd^+$ 
set is a 
set $\CH^+\subseteq P^+_{5,1}$ of oriented exterior lines 
such that 
\begin{itemize}
\item
for every $x\in K$, there are exactly two lines $H_i^+\in \CH^+$, $i=1,2$, such that 
the tangent hyperplane $\pi_5(x)$ contains $\psi H_1^+$ and $\psi H_2^+$, and
\item
the set $\{H_1^+,H_2^+\}$ of these oriented lines depends 
continuously on the point $x$. 
\end{itemize}

b) If $\CH^+$ is an $\it hfd^+$ set, we denote by 
   $$\Pi^+(\CH^+)=\CS^+(\pi_5^+\CH^+)$$ 
the set of all oriented spreads $\CS^+(\pi_5^+ H^+)\subseteq h(K^+)$, $H^+\in \CH^+$, as in Theorem \ref{orsprbyinters}.   
\end{Definition}

In contrast to the non-oriented case, the definition of $\it hfd^+$ sets is quite useless 
without the  condition on continuity of the inverse. This is because it works properly only in 
connection with Lemma \ref{2-1}; compare also Example \ref{patholog}. 
As a compensation, compactness can be deduced 
from this condition. 

\bprop\label{hfdcomp}
Let $\CH^+$ be an $\it hfd^+$ set. Then

a) $\CH^+$  is compact.

b) For each $H^+ \in \CH^+$, there is a point $x \in K$ such that $\psi H^+ \subseteq \pi_5(x)$.
\eprop

\bpf 
For assertion (b), note  that exterior lines are of type $(2,0)$ or $(0,2)$. 
For every $x\in K$, the tangent hyperplane $\pi_5(x)$ contains subspaces of both those types, 
and  the group $\Sigma$ is transitive on the subspaces of either type, whence (b) follows. 

Now assertion (a) follows, because $K$ is compact. Indeed, from the continuity property of 
$\it hfd^+$ sets, we infer that the set of non-ordered pairs $\{H_1^+,H_2^+\}$ of 
oriented lines contained in some tangent hyperplane $\pi_5(x)$ is compact. By (b) it 
follows easily that $\CH^+$ is compact, as well.
\epf

Here is our first main result. It describes \it all \rm topological oriented 
regular parallelisms, whereas the final results of the next section only deal 
with the case that $\dim \mathop \mathrm{span} \CH^+=3$.

\bthm\label{hfd}
If $\CH^+$ is an  $\it hfd^+$ set, then the set 
$\Pi^+(\CH^+) = \CS^+(\pi_5^+\CH^+)$ of oriented 
spreads is a topological oriented regular parallelism, and every topological oriented 
regular parallelism arises in this way.
\ethm

\bcor
Every $\it hfd^+$ set $\CH^+$ is homeomorphic to the 2-sphere $\BS_2$, and it 
consists entirely either of lines of type $(2,0)$ or of lines of type $(0,2)$.
\ecor

\it Proof of Corollary. \rm The polarity $\pi_5^+$ is continuous, and Theorem \ref{orsprbyinters}
asserts that the map $P^+ \to \CS^+(P^+)$ is continuous with respect to the Hausdorff metric. 
Since $\CH^+$ is compact by Proposition \ref{hfdcomp}, it follows that the oriented parallelism 
$\Pi^+(\CH^+)$ is homeomorphic to $\CH^+$, 
and we know that oriented parallelisms are homeomorphic to the 2-sphere. In particular,
$\CH^+$ is connected, and the second assertion follows. \ok
\\

\it Proof of Theorem \ref{hfd}. \rm Using Theorem \ref{orsprbyinters} and 
Proposition \ref{hfdcomp}, we obtain that $\Pi^+ = \Pi^+(\CH^+)$ is a set
of oriented spreads, and that this set is compact with respect to the Hausdorff metric. 
We want to use 
Lemma \ref{2-1} in order to show that every oriented line belongs to exactly 
one of these oriented spreads. At first sight, there seems to be no mapping 
available to which the Lemma 
might be applied. 

However, instead of the condition just stated, is suffices to consider the star
 $\mathfrak{L}_p^+$ of all oriented lines passing through some point $p$ of 
$\Pd$ and to prove that every oriented line in this star
belongs to exactly one spread in $\Pi^+$. We view the star as a subset of $K^+$. 
Now we have the two-sheeted covering map 
$\psi: \mathfrak{L}^+_p \to \mathfrak{L}_p$, and we have a map $\tilde g: \CH^+ \to 
\mathfrak{L}^+_p$ that sends an oriented line $H^+ \in \CH^+$ to the unique oriented 
line of the oriented spread
$\CS^+(\pi_5^+H^+)$ containing $p$.
This map is continuous
by an argument similar to the proof of Theorem \ref{toppar}.
By the properties of the $\it hfd^+$ set and those of the polarity, 
every line  $L \in \mathfrak{L}_x$ belongs to the spreads 
$\psi \CS^+(\pi_5^+H_i^+)$ for exactly two lines $H_1^+, H_2^+ \in \CH^+$ 
and, moreover, the set $\{ H_1^+, H_2^+\}$ depends continuously on $L$.
This means that the composite map $g = \psi \circ \tilde g$ 
has inverse images $g^{-1}(L)$ of cardinality 2, which depend continuously on $L$. 
Now the Lemma tells us that $\tilde g$ is bijective,
which completes the proof that $\Pi^+$ is a compact oriented parallelism. 

It remains to prove the converse, i.e., that every compact oriented parallelism comes from 
some $\it hfd^+$ set. This is obtained without difficulty by retracing all the steps. 
One special point is to show that the set of 
oriented lines from the $\it hfd^+$  set that are contained in a tangent hyperplane 
$\pi_5(x)$, $x \in K$, really depends on $x$ continuously. The reason for 
this is the fact that the 2-valued inverse
map of $\psi: K^+ \to K$ is continuous because $\psi$ is a covering map. Now the 
continuity property of the given oriented parallelism implies that the 
set of two oriented spreads containing the elements of $\psi^{-1}(x)$  depends continuously on
$x\in K$, and this translates to the corresponding property of $\CH^+$.
\ok\\

Every non-oriented $\it hfd$ set yields an $\it hfd^+$ set by taking its inverse 
image with respect to $\psi$, the forgetful map.  $\it hfd^+$ sets obtained in this way will 
be called \it foldable. \rm Clearly we have the following.

\bprop
The oriented parallelism $\Pi^+(\CH^+)$ associated with an $\it hfd^+$ set $\CH^+$ is 
foldable if and only if $\CH^+$ is foldable. \ok
\eprop

\section{Oriented generalized line stars and second main result}\label{GL}

\begin{Definition} \rm  Let $\Pi^* = \Pi^*(\CH^*)$ be an oriented or non-oriented 
compact regular parallelism.

a) The space 
    $$R = R(\Pi^*) := \mathop {\rm span} \CH^*$$
will be called the \it ruler \rm of the parallelism $\Pi^*$.    

b) The number 
$$\dim \Pi^* := \dim R = \dim \mathop {\rm span} \CH^*$$ 
will be called the \it dimension \rm of $\Pi^*$. This terminology was introduced, in the 
non-oriented case, by Betten and Riesinger \cite{hyper}. 
\end{Definition}

If $\dim \Pi^+ = 2$, then $\CH^+ \approx \BS_2$ must consist of all oriented lines in $R$. 
According to \cite{hyper}, Lemma 2.7, $\psi (\Pi^+)$ is the ordinary Clifford parallelism.
Hence $\Pi^+$ is its oriented unfolding, the oriented Clifford parallelism.

In this section, we shall deal with the 
case $\dim \Pi^* = 3$. In this case, passing to a dual object one obtains a 
very convenient description for $\CH^*$. In the non-oriented case, this was 
shown by Betten and Riesinger \cite{Thaswalker}; see \cite{gldirect} for a simple proof. 
We shall see that this direct 
proof carries over to the oriented case almost \it verbatim. \rm The only problem is to capture
the right kind of dual object by a suitable definition. \\

Let an  $\it hfd^+$ set $\CH^+$ with 3-dimensional span be given. 
We know that either all 
lines of $\CH^+$ are
(oriented) $(2,0)$-spaces, or all these lines are $(0,2)$-spaces. 
The two cases are interchanged by 
the map $(x_1,x_2,x_3,y_1,y_2,y_3) \to (y_1,y_2,y_3,x_1,x_2,x_3)$, so without 
loss of generality we may assume that the elements of $\CH^+$ are of type $(0,2)$.
The following result is essentially contained in the proof of \cite{Stevin}, Theorem 23.
We give an abbreviated proof for the sake of completeness. 
 
\bprop
Let $\Pi^* = \Pi^*(\CH^*)$ be a 3-dimensional compact regular parallelism, and assume 
that the elements of $\CH^*$ are of type $(0,2)$ (i.e., negative definite). Then the 
ruler $R(\Pi^*) \in P_{5,3}$ is of type $(1,3)$. In particular, it meets the Klein quadric 
$K$ in an elliptic quadric $Q = K\cap R$.
\eprop

\bpf
If $\CH^*$ were a line star $\mathfrak{L}_y^*$ of $R$, then the point $y$ would belong to all tangent 
hyperplanes $\pi_5(x)$, $x \in K$, a contradiction. Therefore, $R$ is spanned by two lines
$H_1$, $H_2$ from $\psi \CH^*$. Thus $\pi_5R = \pi_5H_1 \wedge \pi_5H_2$ meets $K$ in 
the (empty) intersection of two spreads from $\Pi^*$ and is an exterior line. This implies 
our claim.
\epf

Now we want to start with a potential ruler $R$ of type $(1,3)$ in $\Pf$
and to find an $\it hfd^*$ set $\CH^*$ in $R$. We use the polarity $\pi_3$ of $R$ induced 
by $\pi_5$, in other words, the polarity defined by the restriction to $R$ of the form $f$.
We say that a point of $R$ is \it non-interior \rm with respect to $Q = R\cap K$ 
if it either belongs to $Q$ 
or is contained in a line that misses $Q$. Let $\mathop{Ni}Q$ be the set of 
non-interior points. By a \it 2-secant \rm of $Q$ we mean a line meeting $Q$ in two distinct points.

In the non-oriented case, it is known from \cite{Thaswalker} and \cite{gldirect} that 
the compact $\it hfd$ sets in $R$ are precisely the sets $\pi_3G$, where $G$ is a 
compact set of 2-secants of $Q$ such that every point of $\mathop{Ni}Q$ belongs to 
exactly one line from $G$. A set $G$ of this kind is called a \it generalized line star\rm,
abbreviated \it gl \rm star.
The simplest example is an ordinary line star $G=\mathfrak{L}_x$; then $\CH = \pi_3 \mathfrak{L}_x$ is the 
line set of the plane $\pi_3x$, and we are in the Clifford case. We ask how to define
the correct oriented analogue of a \it gl \rm star. The answer is the following.

\begin{Definition}
\rm Let $Q$ be an elliptic quadric in a 3-dimensional real projective space. 
A set $G^+$
of oriented 2-secants of $Q$ is called an \it oriented gl star \rm or just a
$\it gl^+$ star if every point $x$ of 
$\mathop{Ni}Q$ is incident with exactly two oriented lines from $G^+$, and if the set of 
these two oriented lines depends continuously on the point $x$. 
\end{Definition}

As with $\it hfd^+$ sets, there is no version of the notion of $\it gl^+$ stars 
without the continuity condition, but compactness can be deduced from it.

\ble
Every $\it gl^+$ star $G^+$ is compact, but compactness cannot replace the continuity 
condition in the above definition.
\ele

\bpf
Compactness of $G^+$ follows directly from compactness of $Q$ via the continuity property.
That compactness does not conversely imply the continuity condition will be demonstrated 
by Example \ref{patholog}.
\epf

\bthm\label{gl}
Let $R$ be a 3-space of type $(1,3)$ in $\Pf$. The $\it hfd^+$ sets in 
$R$ are precisely the sets $\CH^+ = \pi_3^+G^+$, where $G^+$ is a $\it gl^+$ star 
with respect to $Q = R\cap K$. The space $R$ is generated by $\CH^+$ if and only if 
$G^+$ is not an ordinary star of oriented lines.
\ethm

Combining this with Theorem \ref{hfd} we obtain our main result for the 3-dimensional case:

\bcor\label{main}
The compact oriented regular parallelisms of $\Pd$ of dimension $d \le 3$ are 
precisely the parallelisms 
   $$\Pi^+(G^+) := \Pi^+(\pi_3^+G^+),$$
where $G^+$ is a $\it gl^+$ star in a 3-space $R$ of $\Pf$ that meets the Klein 
quadric $K$ in an elliptic quadric of $R$.  The parallelism is Clifford if and only if 
$G^+$ is an ordinary star of oriented lines. \ok
\ecor

If the ruler $R$ is of type $(1,3)$, as we have assumed previously, then 
the 3-spaces $\pi_5^+\pi_3^+L^+$, $L^+ \in G^+$, which define the oriented spreads of
$\Pi^+(G^+)$, are of type $(3,1)$. \\

The proof of Theorem \ref{gl} uses the following Lemma.

\ble\label{conic}
\rm (\cite{gldirect}, 2.4) \it  Let $U \in P_{3,1}$ be a subspace of type $(1,3)$ and 
consider the elliptic 
quadric $Q = X \cap K$. For every $x\in K\setminus Q$, the tangent hyperplane $\pi_5(x)$
intersects $Q$ in a non-degenerate conic, and every non-degenerate conic in $Q$ arises 
in this way. 
\ele

\bpf
The quadric $Q$ represents a spread in $\Pd$, hence the line represented by 
$x \notin Q$ intersects infinitely many lines in this spread. This means that $x$ is 
$f$-orthogonal to infinitely many elements $q \in Q$. These elements then belong to 
the plane $\pi_5(x)\wedge R$.
The converse follows using transitivity properties of the group $\Sigma$.
\epf

\it Proof of Theorem \ref{gl}. \rm The proof is practically the same as the proof of 
\cite{gldirect}, Theorem 2.3. Only the words `exactly one' have to be replaced by `exactly two' 
where appropriate. For the sake of completeness, we give the proof that a $\it gl^+$ star 
$G^+$ yields an $\it hfd^+$ set; the converse direction is similar. 

For $L^+ \in G^+$ and $x \in K$, we have that 
$\psi H^+ = \psi \pi_3^+(L^+)\in \psi \CH^+$ is contained in the 
tangent hyperplane $\pi_5(x)$ if and only if it is contained in $\pi_5(x) \cap R$, and 
this happens if and only if $\psi L^+$ contains the point $r(x) := \pi_3(\pi_5(x) \cap R)$.
If $x\notin Q$, then $r(x) \in \mathop{Ni}Q$ by Lemma \ref{conic}. If $x \in Q$, then again,
$r(x) = x \in \mathop{Ni}Q$. There are exactly two lines 
$L_i^+ \in G^+$, $i= 1,2$, containing $r(x)$. The set $\{L_1^+, L_2^+\}$ depends continuously
on $r(x)$, which in turn is a continuous function of $x$.  \ok\\

Every compact non-oriented $\it gl$ star $G$ has the continuity property, that is, 
the line from $G$ passing through a non-interior point $p$ continuously depends on $p$, 
see \cite{gldirect}, Theorem 3.2.
This implies that taking the inverse image of  $G$ with respect to $\psi$, 
we obtain a $\it gl^+$ star $G^+$. Such examples are said to be \it foldable. \rm Clearly, 
we have the following.

\bprop
Let $G^+$ be a $\it gl^+$ star. Then $G^+$ and the associated $\it hfd^+$ set 
$\pi_3^+G^+$ as well as the associated oriented parallelism $\Pi^+(G^+)$ are either 
all foldable or all non-foldable.
\eprop

\section{Examples}\label{ex}

With a non-oriented \it gl \rm star $G$, there is associated the involutory 
homeomorphism $\sigma : Q \to Q$ which sends a point $x \in Q$ to the second 
point of intersection of the line $G_x \in G$ containing $x$. This involution 
carries all information about $G$. In the oriented case, the intersection points 
of $L^+ \in G^+$ with $Q$ can be distinguished: at one point, $L^+$ enters the closed 
ball $B$ bounded by $ Q$ (i.e., a positive tangent vector points inward), and at the 
other point, the line leaves $B$. Let us call these points $e(G^+)$
and $l(G^+)$, respectively. 

\ble
If $G^+$ is a $\it gl^+$ star with respect to $Q$, then at every point of $Q$ exactly one 
oriented line $L^+ \in G^+$ enters the closed ball $B$ bounded by $Q$, and exactly one 
leaves $B$.
\ele

\bpf
As every oriented line $L^+ \in G^+$ has an entry point and a leave point, 
both the set of entry points and the set of leave points are nonempty. Lines entering 
at $x_n \to x$ cannot converge to a line leaving 
at $x$, hence both sets are closed. So the connected quadric $Q$ 
is a disjoint union of three closed sets: the set $EL$ of points where one line enters and 
one line leaves, the set $EE$ where two lines enter, and the set $LL$ where two lines leave.
That these sets are closed follows from the continuity property in the 
definition of $\it gl^+$ stars. Only one of the three sets can be nonempty, and 
this set can only be $EL$. 
\epf

\begin{Definition} \rm
Given a $\it gl^+$ star $G^+$, define a map $\rho: Q \to Q$ by sending $x\in Q$ to the 
leave point of the unique line $G_x^+$ that enters the ball $B$ at $x$.
We call this map the \it characteristic map \rm of $G^+$.
\end{Definition}

The following lemma is now obvious.

\ble\label{charmap}
If $G^+$ is a $\it gl^+$ star with characteristic map $\rho$, then $\rho$ is a 
fixed point free homeomorphism, 
and $G^+$ is the set of all lines
$G_x^+ = x \vee \rho(x)$, oriented in such a way that the interval 
$G_x \cap B$ is traversed from $x$ to $\rho(x)$. \ok
\ele

This opens up a huge set of candidates for $\it gl^+$ stars. Every fixed 
point free homeomorphism  of the 
2-sphere may be tested for its potential of defining a $\it gl^+$ star. 
The two oriented lines passing through $x \in Q$ are 
then $x \vee \rho(x)$ and $\rho^{-1}(x) \vee x$, so the continuity condition for the 
set of oriented lines of the $\it gl^+$ star passing through a given point is 
satisfied on $Q$ at least. 
In general, the test will fail nevertheless, 
but one may suspect that there are far more successes than with ordinary $gl$ stars,
where the characteristic map is an involution. We pursue this a bit further.

\begin{Definition} \rm
An (oriented) $\it gl^+$ star is said to be \it foldable \rm if by forgetting orientations 
it yields an ordinary $\it gl$ star. As before, we call an oriented 
parallelism $\Pi$ \it foldable \rm if by forgetting orientations we get an ordinary 
parallelism.
\end{Definition}

Again we have two obvious facts:

\bprop\label{fold}
a) A $\it gl^+$ star is foldable if and only if its characteristic map is an involution.

b)
The oriented regular parallelism defined by a $\it gl^+$ star is foldable if and only if the 
$\it gl^+$ star is foldable.
\ok
\eprop

\begin{Example}\label{nonfold}
\rm Here is a class of non-foldable $\it gl^+$ stars. Compare also Example \ref{patholog}.
In \cite{torus},
a large set of rotationally symmetric $\it gl$ stars were constructed. They are defined by
their characteristic involutions $\sigma: \BS_2 \to \BS_2$. On the equator ($z=0$) the involutions
agree with the antipodal map, that is, $\sigma (x,y,0) = - (x,y,0)$. Now we take two 
different involutions $\sigma_1$ and $\sigma_2$ of this kind and define a homeomorphism
$\rho \BS_2 \to \BS_2$ by sending $(x,y,z)$ to $\sigma_1(x,y,z)$ if $z \ge 0$, and to
$\sigma_2(x,y,z)$ if $z \le 0$. Then it is easily checked that this map $\rho$ is 
the characteristic map of a non-foldable parallelism.
\end{Example}

\bthm 
The examples described above define compact oriented regular parallelisms. 
These parallelisms  are 
non-foldable, and they admit a 2-dimensional torus group of automorphisms.
\ethm

The \it proof \rm \, is almost automatic. 
For the automorphism group, 
compare \cite{Stevin}, Theorem 31 or 
\cite{torus}, Proposition 4.1. We note that a 2-torus is as much symmetry as an 
oriented regular parallelism can have without being Clifford, see \cite{torus}, Theorem 2.1.
\\

If a $\it gl^+$ star is to be constructed from its characteristic map, then the 
defining incidence condition can be relaxed and the continuity property can be replaced by 
an orientation rule. The analogous result for ordinary $\it gl$ stars is Proposition 5.1 of 
\cite{torus}. 

\bthm
Let $Q$ be an elliptic quadric in a real projective 3-space $R$ and let $\rho: Q\to Q$ 
be a fixed point free homeomorphism. Let $G^+=G^+(\rho)$ be the set of all oriented lines
$G_x^+ = x \vee \rho(x)$, oriented in such a way that the interval 
$\psi(G_x^+) \cap B$ is traversed from $x$ to $\rho(x)$.

If each point of $\mathop{Ni}Q$ 
is incident with at most two of these oriented lines, then $G^+$ is a $\it gl^+$ star. 
In particular,  
the continuity property of $\it gl^+$ stars comes for free in this situation.
\ethm

\bpf
1. We may assume that $Q$ is the unit sphere in the affine space $\BR^3$, of which $R$ 
is the projective closure. By construction, $G^+$ is compact.
First we look at `affine' points $a \in A :=\BR^3 \cap \mathop{Ni} Q$.
The orientation of a line $L^+ \in G^+$ defines an order relation on the affine line
$\psi L^+ \cap \BR^3$. By the construction of $G^+$, the entry and leave points of $L^+$
satisfy $e(L^+) < l(L^+)$. We call $L^+$ a positive line with respect to $a\in L^+$ if 
$l(L^+) \le a$.  The only other possibility is $a \le e(L^+)$, in which case
we call $L^+$ a negative line with respect to $a$. 

2. Let $L_n^+\in G^+$ be  a sequence of oriented
lines that are positive with respect to points $a_n$. If both sequences converge, then
$\lim L_n^+$ is positive with respect to $\lim a_n$. Let $B$ be the set of points incident
both with a positive line and with a negative line from $G^+$. If $a_n\in B$, $a_n\to a$ 
and if the lines $L_n^+ \in G^+$ are positive with respect to $a_n$, then these 
lines accumulate at some positive $L^+$ line for $a$, by compactness of $G^+$, and 
similarly for negative lines. Hence, we have
$a \in B$. It follows that in fact
$L_n^+$ converges to $L^+$, 
and the desired continuity condition is satisfied on $B$.
We proceed to show that $B = A$.

3. It suffices to show that every point of $A$ is on a negative line. For $x \in Q$, 
we have the oriented line $G_x^+ \in G^+$.  We
define a ray $W_x$ as the connected component of $\psi (G_x^+)\cap A$ that contains $x$.
For $r \ge 1$, let $g_r(x)$ be the intersection point of $W_x$ with the sphere $rQ$ 
of radius $r$. Then $G_x^+$ is a negative line with respect to $g_r(x)$. Let $h :A\to Q$
be the map sending $y$ to $y\Vert y\Vert^{-1}$. Then $k_r := h\circ g_r$, $r \ge 1$, 
is a family of
pairwise homotopic maps $Q\to Q$, with $k_1 = \rm id$. So all of these maps have mapping degree one
and, hence, are surjective. This proves our claim, and the continuity condition is proved
on the set $A$.

4. It remains to prove continuity for points on the plane $I$ at infinity. Every point
$y \in I$ has a neighborhood $U$ homeomorphic to $\BR^2 \times [-1,1]$ such that 
$U\cap I$ corresponds to $\BR^2 \times \{0\}$ and that $U \cap rQ$ for some large number $r$
corresponds to $\BR^2 \times \{-1,1\}$. By compactness of $Q$, and for large enough $r$,
there is a neighborhood 
$V\subseteq U$ of $y$ such that each oriented line meeting 
both $V$ and $Q$
intersects both the top layer $\BR^2 \times \{1\}$ and the bottom 
layer $\BR^2 \times \{-1\}$ of $U$.
These oriented lines may therefore be divided into the set of upward lines, 
traversing $U$ from bottom to top, and downward lines. By the previous steps, each point
of $V \setminus I$ is incident with both an upward line 
and a downward line from $G^+$.
As before, we may conclude that the same is true for all points of $V$, and the 
continuity condition is guaranteed for these points. 
\epf

We conclude with an example demonstrating the necessity of the 
condition on continuity of the inverse in the definitions of $\it gl^+$ stars 
and, hence, also in that of $\it hfd^+$ sets. It shows that compactness is not a possible 
surrogate condition.

\begin{Example} \label{patholog}
\rm We start by defining an ordinary $\it gl$ star $G_1$ with respect to $Q = \BS_2$ in $\Pd$. 
It contains all lines of the plane $z=0$ that pass through the origin $o = (0,0,0) \in \BR^3$, 
but no other lines  containing $o$. This $\it gl$ 
star does not have rotational symmetry. It has to be so, because rotationally
symmetric $\it gl$ stars inevitably contain the rotation axis. Incidentally, this 
is the simplest known example with this property. For more such examples, see Section 7.2 of 
\cite{aut}.

Consider the points 
    $$p_t = (\sqrt{1-t^2}, 0, t) \in \BS_2 \quad and \quad q_t = (f(t),0,0),$$
where $t \in [0,\frac{1}{2}]$ and where $f: [0,\frac{1}{2}]\to [\frac{1}{2} ,0]$ is a 
strictly decreasing bijection.  Let $L_t$ be the line $p_t \vee q_t$. Then the slope of $L_t$ 
strictly decreases from 0 to $-\infty$, and $L_{\frac{1}{2}}$ is parallel to the 
$z$-axis $Z$. Now rotate the line $L_t$ about the axis $A_t$ parallel to $Z$ and 
passing through $q_t$. These rotated lines fill a cone $C_t$, and  we define $G_1$ 
to be the set of all
lines obtained by rotating $L_t$ for all $t\in [0,\frac{1}{2}]$. 
The non-interior part $C_t\cap \mathop {\it Ni} Q$ lies completely 
inside $C_s$ for all $s < t$. 
By continuity, every point of $ \mathop {\it Ni} Q$ is incident with exactly one line from 
$G_1$, and we have defined a $\it gl$ star.    

Now consider the ordinary star of lines $G_2 = \frak L_o$, and observe that $G_1\cap G_2$ 
consists precisely of the horizontal lines passing through $o$. 
This gives us two possibilities to try forming a $\it gl^+$ star. 

1. We take the horizontal lines with both orientations, the lines of $G_1$ with 
upward orientation, and the lines of $G_2$ with downward orientation. 
Like Example \ref{nonfold}, this gives 
a nice non-foldable $\it gl^+$ star satisfying the continuity condition, 
but without rotational symmetry.

2. We take the elements of $G_1 \cup G_2$, all of them with upward orientation 
(and two orientations for the horizontal ones). This is a compact set 
$A\approx \BS_2$ of oriented lines such that
every point $p \in \mathop{Ni} Q$ lies on precisely two 
distinct oriented lines in $A$ (here we use the
information about the intersection $G_1\cap G_2$). But the 
set of those two lines does not depend on $p$ continuously when $p$ is in the plane $z = 0$.
So $A$ is \it not \rm a $\it gl^+$ star. By Corollary \ref{main}, it does  not 
define an oriented parallelism.
\end{Example}

The possibility of such examples quite puzzled the author until he understood 
the relevance of the continuity condition. What the example tells us is that, 
in contrast with non-oriented
$\it gl$ stars, compactness does not suffice to ensure the continuity property of 
a $\it gl^+$ star. Yet for oriented parallelisms themselves, compactness is 
enough, by Proposition \ref{toppar}.

\bibliographystyle{plain}

\noindent{Rainer L\"owen\\ 
Institut f\"ur Analysis und Algebra\\
Technische Universit\"at Braunschweig\\
Universit\"atsplatz 2\\
38106 Braunschweig\\
Germany

\end{document}